# On the Relation of IOS-Gains and Asymptotic Gains For Linear Systems

## Iasson Karafyllis

[*]Dept. of Mathematics, National Technical University of Athens,
Zografou Campus, 15780, Athens, Greece.
emails: iasonkar@central.ntua.gr , iasonkaraf@gmail.com

**Abstract**

This paper presents a fundamental relation between Output Asymptotic Gains (OAG) and Input-to-Output Stability (IOS) gains for linear systems. For any Input-to-State Stable, strictly causal linear system the minimum OAG is equal to the minimum IOS-gain. Moreover, both quantities can be computed by solving a specific optimal control problem and by considering only periodic inputs. The result is valid for wide classes of linear systems (involving delay systems or systems described by PDEs). The characterization of the minimum IOS-gain is important because it allows the non-conservative computation of the IOS-gains, which can be used in a small-gain analysis. The paper also presents a number of cases for finite-dimensional linear systems, where exact computation of the minimum IOS gain can be performed.

## 1. Introduction

Input-to-State Stability (ISS) has played a crucial role for the development of modern nonlinear control theory. Since its first formulation for finite-dimensional systems by E. D. Sontag in [16], ISS has been extended to cover many important cases. The notion of Input-to-Output Stability (IOS) for finite-dimensional systems has been studied in [1,5,20]. Extensions of ISS/IOS to delay systems and PDEs were recently studied in [5,6,7,10,11,12,14,15] (but see also references therein).

The study of ISS/IOS involves two important notions that play a significant role in the description of the properties of a control system: the notion of the gain of an input and the notion of the asymptotic gain of an input. The notion of the gain of an input has been used extensively for the derivation of stability properties by means of the Small-Gain Theorem (see [3,5]). Asymptotic gain properties for finite-dimensional systems were introduced in [1,2,10,19]. More specifically, in [19] it was shown the ISS superposition theorem, which was extended in [1] for the case of the output asymptotic gain property. Recently, the asymptotic gain property has been used in time-delay systems (see [6]), systems described by Partial Differential Equations (PDEs) (see [8,9]) and abstract infinite-dimensional systems (see [10,12]). For linear systems, where the asymptotic gain function and the gain function are linear functions, the minimum (asymptotic) gain (coefficient) may be used as a measure of the sensitivity of the system with respect to external disturbances.

This paper presents a fundamental relation between output asymptotic gains and IOS-gains for linear systems. As mentioned above, for linear systems it is possible to define the minimum IOS-gain and the minimum Output Asymptotic Gain (OAG) of an input: both are well-defined quantities, which can be defined independently of each other. The present paper offers a fundamental relation for ISS, strictly causal linear systems: the minimum OAG is equal to the minimum IOS-gain. Moreover, both quantities can be computed by solving a specific optimal control problem and by considering only periodic inputs. The characterization of the minimum IOS-gain is important because it allows the non-conservative computation of the IOS-gains, which can be used in a small-gain analysis. Indeed, the present paper offers a number of cases for finite-dimensional linear systems, where exact computation of the minimum IOS-gain can be performed. However, the fundamental result about the equality of the minimum IOS-gain and the minimum OAG of an input is valid for wide classes of linear systems (involving delay systems or systems described by PDEs).

The structure of the paper is as follows. All fundamental notions for linear control systems are provided in Section 2. Section 3 of the paper contains a motivating example of a linear delay system, where a straightforward analysis indicates that the minimum OAG may be strictly smaller than the minimum IOS gain. The main results of the paper, which are stated and proved in Section 4, show that this cannot be the case and consequently, such an indication for the delay system presented in Section 3 is only an artifact of the stability analysis. Section 5 of the paper is devoted to cases of finite-dimensional linear systems, where exact computation of the minimum IOS-gain can be performed or useful estimates of the IOS-gain can be obtained. The concluding remarks of the present work are provided in Section 6.

**Notation.** Throughout this paper, we adopt the following notation.

* $\Re_+ := [0, +\infty)$. For $x \in \Re$, $[x]$ denotes the integer part of $x$, i.e., the largest integer which is less or equal to $x$.

* Let $S \subseteq \Re^n$ be an open set and let $A \subseteq \Re^n$ be a set that satisfies $S \subseteq A \subseteq cl(S)$. By $C^0(A; \Omega)$, we denote the class of continuous functions on $A$, which take values in $\Omega \subseteq \Re^m$. By $C^k(A; \Omega)$, where $k \geq 1$ is an integer, we denote the class of functions on $A \subseteq \Re^n$, which takes values in $\Omega \subseteq \Re^m$ and has continuous derivatives of order $k$. In other words, the functions of class $C^k(A; \Omega)$ are the functions which have continuous derivatives of order $k$ in $S = \text{int}(A)$ that can be continued continuously to all points in $\partial S \cap A$. When $\Omega = \Re$ then we write $C^0(A)$ or $C^k(A)$.

* Let $a < b$ be given real numbers. For every $p \geq 1$, $L^p(a,b)$ denotes the equivalence class of Lebesgue measurable functions $f:[a,b] \to \Re$ for which $\|f\|_p = \left( \int_a^b |f(x)|^p \, dx \right)^{1/p} < +\infty$. $L^\infty(a,b)$ denotes the equivalence class of measurable functions $f:[a,b] \to \Re$ for which $\|f\|_\infty := \operatorname*{ess\,sup}_{x \in (a,b)} (|f(x)|) < +\infty$.

* Let $u: \Re_+ \times [0,1] \to \Re$ be given. We use the notation $u[t]$ to denote the profile at certain $t \geq 0$, i.e., $(u[t])(x) = u(t,x)$ for all $x \in [0,1]$.

* For a vector $x \in \Re^n$ we denote by $|x|$ its usual Euclidean norm and by $x'$ its transpose. By $|A| := \sup\{|Ax|; x \in \Re^n, |x| = 1\}$ we denote the induced norm of a matrix $A \in \Re^{m \times n}$ and $I_n \in \Re^{n \times n}$ denotes the $n \times n$ identity matrix, while $A' \in \Re^{n \times m}$ denotes the transpose of $A \in \Re^{m \times n}$. A square matrix $A \in \Re^{n \times n}$ is called a Metzler matrix if all its off-diagonal entries are nonnegative.

* For an integer $k \geq 1$, $H^k(0,1)$ denotes the Sobolev space of functions in $L^2(0,1)$ with all its weak derivatives up to order $k \geq 1$ in $L^2(0,1)$.

* Let $D \subseteq \Re^l$ be a non-empty set and $A \subseteq \Re_+$ an interval. By $L^\infty_{loc}(A; D)$ we denote the class of all Lebesgue measurable and locally bounded mappings $d: A \to D$.

## 2. Notions for Linear Systems

In this work we study control systems which satisfy the requirements of the following definition.

**Definition 2.1 (Control System):** *A control system is a quintuplet $\Sigma = (X, Y, M(U), \phi, h)$ consisting of*

*(i) a normed vector space $(X, \|\ \|_X)$, called the state space, endowed with the norm $\|\ \|_X$,*

*(ii) a normed vector space $(Y, \|\ \|_Y)$, called the output space, endowed with the norm $\|\ \|_Y$,*

*(iii) a normed vector space $(U, \|\ \|_U)$, called the space of input values, endowed with the norm $\|\ \|_U$,*

*(iv) a set $M(U)$ of locally bounded functions $u: \Re_+ \to U$, which satisfies the following requirement: for all $u_1, u_2 \in M(U)$ and all $\tau > 0$, the functions $\tilde{u}: \Re_+ \to U$, $\hat{u}: \Re_+ \to U$, $\delta_\tau u_1 : \Re_+ \to U$ defined by $\tilde{u}(t) := u_1(t)$ for $t \in [0, \tau]$, $\tilde{u}(t) := u_2(t - \tau)$ for $t > \tau$, $\hat{u}(t) := u_1(t)$ for $t \in [0, \tau)$, $\hat{u}(t) := u_2(t - \tau)$ for $t \geq \tau$ and $(\delta_\tau u_1)(t) := u_1(t + \tau)$ for $t \geq 0$ belong to $M(U)$,*

*(v) a map $h: X \to Y$, called the output map,*

*(iv) a map $\phi: D_\phi \to X$, where $D_\phi \subseteq \Re_+ \times X \times M(U)$ (called the transition map), which satisfies the following properties:*

   *(a) (existence) for all $(x, u) \in X \times M(U)$ there exists $t > 0$ such that $[0, t] \times \{(x, u)\} \subseteq D_\phi$,*

   *(b) (identity) for all $(x, u) \in X \times M(U)$ it holds that $\phi(0, x, u) = x$,*

   *(c) (causality) for every $(t, x, u) \in D_\phi$ and for every $\tilde{u} \in M(U)$ with $\tilde{u}(s) = u(s)$ for all $s \in [0, t]$, it holds that $(t, x, \tilde{u}) \in D_\phi$ and $\phi(t, x, u) = \phi(t, x, \tilde{u})$,*

   *(d) (semigroup) for every $(t, x, u) \in D_\phi$ with $t > 0$ it holds that $(\tau, x, u) \in D_\phi$ for all $\tau \in [0, t]$ and $(t - \tau, \phi(\tau, x, u), \delta_\tau u) \in D_\phi$ with $\phi(t, x, u) = \phi(t - \tau, \phi(\tau, x, u), \delta_\tau u)$.*

Definition 2.1 is similar to other definitions of control systems given in [5,12,17]. However, notice that contrary to [12], we do not assume continuity of the mapping $t \to \phi(t, x, u)$.

In our work, we work with linear control systems, which are control systems that satisfy additional properties.

**Definition 2.2 (Linear Control System):** *A control system $\Sigma = (X, Y, M(U), \phi, h)$ is called a linear control system, if*
*(i) the output map $h: X \to Y$ is continuous and linear,*
*(ii) for all $u_1, u_2 \in M(U)$ and all $\lambda \in \Re$, the functions $\hat{u}: \Re_+ \to U$, $\breve{u}: \Re_+ \to U$ defined by $\hat{u}(t) := u_1(t) + u_2(t)$ and $\breve{u}(t) := \lambda u_1(t)$ for $t \geq 0$ belong to $M(U)$ and are denoted by $(u_1 + u_2) \in M(U)$ and $(\lambda u_1) \in M(U)$, respectively,*
*(iii) $D_\phi = \Re_+ \times X \times M(U)$,*
*(iv) for all $u_1, u_2 \in M(U)$, $x_1, x_2 \in X$ and all $\lambda, \mu \in \Re$, it holds that $\phi(t, \lambda x_1 + \mu x_2, \lambda u_1 + \mu u_2) = \lambda \phi(t, x_1, u_1) + \mu \phi(t, x_2, u_2)$ for all $t \geq 0$,*
*(v) for all $v \in U$ there exists $u \in M(U)$ with $u(0) = v$. Moreover, for every $u \in M(U)$ and $T > 0$, the $T$-periodic extension of $u$ given by $w(t) = u\left(t - T\left[\dfrac{t}{T}\right]\right)$ for $t \geq 0$ is a function of class $M(U)$.*

Definition 2.2 allows us to consider:

1) Finite-dimensional linear time-invariant systems of the form $\dot{x} = Ax + Bu$, $y = Cx$, where $x \in \Re^n$, $u \in \Re^m$, $y \in \Re^p$ and $A \in \Re^{n \times n}$, $B \in \Re^{n \times m}$, $C \in \Re^{p \times n}$. In this case we have $X = \Re^n$, $Y = \Re^p$ while $U$ can be any subspace of $\Re^m$ and $M(U) = L^\infty_{loc}(\Re_+; U)$ is the set of all Lebesgue measurable and locally essentially bounded inputs $u : \Re_+ \to U$.

2) Time-delay systems of the form $\frac{d}{dt}(x(t) - Qx_{t-\delta}) = Ax_t + Bu(t)$, $y_t = Cx_t$, where $\delta > 0$ is a constant $x(t) \in \Re^n$, $u(t) \in \Re^m$, $(x_t)(\theta) = x(t+\theta)$ for $\theta \in [-\tau, 0]$ with $\tau > 0$, $x_t \in C^0([-\tau, 0]; \Re^n)$, $y_t \in C^0([-r, 0]; \Re^p)$ with $r \in [0, \tau]$ (we are using the convention $C^0([-r, 0]; \Re^p) = \Re^p$ for $r = 0$), $A, Q : C^0([-\tau, 0]; \Re^n) \to \Re^n$, $C : C^0([-\tau, 0]; \Re^n) \to C^0([-r, 0]; \Re^p)$ are bounded linear operators, and $B \in \Re^{n \times m}$ is a real matrix. Notice that this class of delay systems includes conventional delay systems (when $Q \equiv 0$) as well as neutral delay systems in Hale's form (see [15]). In this case we have $X = C^0([-\tau - \delta, 0]; \Re^n)$, $Y = C^0([-r, 0]; \Re^p)$ while $U$ can be any subspace of $\Re^m$ and $M(U) = L^\infty_{loc}(\Re_+; U)$ is the set of all Lebesgue measurable and locally essentially bounded inputs $u : \Re_+ \to U$ (see [15]).

3) Integral delay equations of the form (see [6]) $x(t) = A\breve{x}_t + Bu_t$, $y(t) = C\breve{x}_t$, where $x(t) \in \Re^n$, $u(t) \in U \subseteq \Re^m$, $y(t) \in \Re^p$, $r > 0$, $U$ is a subspace of $\Re^m$, $\breve{x}_t \in L^\infty([-r, 0); \Re^n)$, $u_t \in L^\infty([-r, 0]; U)$ are defined by $(\breve{x}_t)(s) = x(t+s)$ for $s \in [-r, 0)$, $(u_t)(s) = u(t+s)$ for $s \in [-r, 0]$, $A : L^\infty([-r, 0); \Re^n) \to \Re^n$, $B : L^\infty([-r, 0]; \Re^m) \to \Re^n$, $C : L^\infty([-r, 0); \Re^n) \to \Re^p$ are bounded linear operators for which there exist constants $N, M > 0$ such that

$$|A\breve{x}| \leq N\theta \sup_{-\theta \leq s < 0} |\breve{x}(s)| + M \sup_{-r \leq s < -\theta} |\breve{x}(s)|, \text{ for all } \theta \in (0, r) \quad (2.1)$$

and have the property that for every $\delta > 0$, $\breve{\eta} \in L^\infty([-r, \delta); \Re^n)$, $w \in L^\infty([-r, \delta]; U)$, the function $z : [-r, \delta) \to \Re^n$ defined by $z(t) = \breve{\eta}(t)$ for $t \in [-r, 0)$ a.e. and $z(t) = A\breve{\eta}_t + Bu_t$ for $t \in [0, \delta)$ a.e. satisfies $z \in L^\infty([-r, \delta); \Re^n)$. In this case we have $X = L^\infty([-r, 0); \Re^n)$, $Y = \Re^p$ and $M(U) = L^\infty_{loc}(\Re_+; U)$ is the set of all Lebesgue measurable and locally essentially bounded inputs $u : \Re_+ \to U$ (see [6]).

4) Many classes of systems described by PDEs. For example, consider the parabolic PDE problem with a boundary input

$$\frac{\partial w}{\partial t}(t, x) - \frac{1}{r(x)} \frac{\partial}{\partial x}\left(p(x) \frac{\partial w}{\partial x}(t, x)\right) + \frac{q(x)}{r(x)} w(t, x) = 0, \; x \in (0, 1) \quad (2.2)$$

$$b_1 w(t, 0) + b_2 \frac{\partial w}{\partial x}(t, 0) = a_1 w(t, 1) + a_2 \frac{\partial w}{\partial x}(t, 1) - u(t) = 0, \quad (2.3)$$

where $w(t, x) \in \Re$, $u(t) \in \Re$, $p \in C^1([0, 1]; (0, +\infty))$, $r \in C^0([0, 1]; (0, +\infty))$, $q \in C^0([0, 1]; \Re)$, $a_1, a_2, b_1, b_2$ are real constants with $|a_1| + |a_2| > 0$, $|b_1| + |b_2| > 0$. Let $\lambda_n \in \Re$ ($n = 1, 2, ...$) denote the eigenvalues of the Sturm-Liouville operator $A : D \to L^2_r(0, 1)$ defined by

$$(Af)(x) = -\frac{1}{r(x)} \frac{d}{dx}\left(p(x) \frac{df}{dx}(x)\right) + \frac{q(x)}{r(x)} f(x), \text{ for all } f \in D \text{ and } x \in (0, 1) \quad (2.4)$$

where $D \subseteq H^2(0,1)$ is the set of all functions $f:[0,1] \to \Re$ for which $b_1 f(0) + b_2 \frac{df}{dx}(0) = a_1 f(1) + a_2 \frac{df}{dx}(1) = 0$. The eigenvalues of $A$ form an infinite, increasing sequence $\lambda_1 < \lambda_2 < ... < \lambda_n < ...$ with $\lim_{n \to \infty}(\lambda_n) = +\infty$ and to each eigenvalue corresponds exactly one eigenfunction $\phi_n \in C^2([0,1]; \Re)$ that satisfies $A\phi_n = \lambda_n \phi_n$ and $b_1 \phi_n(0) + b_2 \frac{d\phi_n}{dx}(0) = a_1 \phi_n(1) + a_2 \frac{d\phi_n}{dx}(1) = 0$. We assume that $\sum_{n=N}^{\infty} \lambda_n^{-1} \max_{0 \le x \le 1}(|\phi_n(x)|) < +\infty$, for certain $N > 0$ with $\lambda_N > 0$. Let $M(\Re)$ denote the set of locally bounded functions $u: \Re_+ \to \Re$ for which there exists an increasing sequence of times $\{\tau_i \ge 0, i = 0,1,2,...\}$ with $\tau_0 = 0$, $\lim_{i \to +\infty}(\tau_i) = +\infty$ such that: (i) $u \in C^2(I)$, where $I = \Re_+ \setminus \{\tau_i \ge 0, i = 0,1,2,...\}$, (ii) for every $\tau_i > 0$ the left limits of $u(t), \dot{u}(t), \ddot{u}(t)$ when $t$ tends to $\tau_i$ are finite, (iii) for every $\tau_i \ge 0$ the right limit of $u(t)$ when $t$ tends to $\tau_i$ is finite, and (iv) $\sup_{t \in (\tau_i, \tau_{i+1})}(|\dot{u}(t)|) < +\infty$ for $i = 0,1,2,...$. Theorem 4.10 on page 89 in [7] guarantees that for every $w_0 \in L^2(0,1)$ and $u \in M(\Re)$ there exists a unique mapping $w \in C^0(\Re_+; L^2(0,1))$ with $w \in C^1(I \times [0,1])$ satisfying $w[t] \in C^2([0,1])$ for all $t > 0$, $\lim_{t \to \tau_i^-}(w(t,x)) = w(\tau_i, x)$, $\lim_{t \to \tau_i^-}\left(\frac{\partial w}{\partial t}(t,x)\right) = -(Aw[\tau_i])(x)$, $\lim_{t \to \tau_i^-}\left(\frac{\partial w}{\partial x}(t,x)\right) = \frac{\partial w}{\partial x}(\tau_i, x)$ for all $i \ge 1$ and $x \in [0,1]$, $w[0] = w_0$, (2.2) for all $(t,x) \in I \times (0,1)$ and (2.3) for all $t \in I$. Define for every $w_0 \in L^2(0,1)$ and $u \in M(\Re)$ the transition map $\phi(t, w_0, u) = w[t]$. In this case, we have a linear control system with state space $X$ being any normed linear space that satisfies $L^2(0,1) \subseteq X \subseteq C^2([0,1])$ endowed with the $L^2(0,1)$ norm (although different norms may be also considered depending on the choice of $X$; for example, if $X = C^0([0,1])$ then we may use the sup norm). Moreover, in this case we have $U = \Re$, $Y = X$ and $h(x) = x$. Similarly, we can treat parabolic PDE problems with distributed inputs.

For linear systems covered by Definition 2.2, we make use of the stability notions given in the following definition.

**Definition 2.3 (Stability Notions):** **(a)** *A linear control system* $\Sigma = (X, Y, M(U), \phi, h)$ *is Output Exponentially Stable (OES) if there exist constants* $M, \sigma > 0$ *such that for every* $x \in X$ *and* $t \ge 0$ *it holds that*

$$\|h(\phi(t,x,0))\|_Y \le M \exp(-\sigma t) \|x\|_X \tag{2.5}$$

**(b)** *A linear control system* $\Sigma = (X, Y, M(U), \phi, h)$ *is called Exponentially Stable (ES) if it is OES with output the identity map* $h(x) = x$.

**(c)** *A linear control system* $\Sigma = (X, Y, M(U), \phi, h)$ *is called Exponentially Input-to-Output Stable (exp-IOS), if there exist constants* $M, \sigma > 0$, $g \ge 0$ *such that for every* $(x, u) \in X \times M(U)$ *and* $t \ge 0$ *it holds that*

$$\|h(\phi(t,x,u))\|_Y \le M \exp(-\sigma t) \|x\|_X + g \sup_{0 \le s \le t}\left(\|u(s)\|_U\right) \tag{2.6}$$

*The constant* $g \geq 0$ *is called an IOS-gain of the input* $u$. *If* $U \neq \{0\}$ *then the constant* $\gamma \geq 0$ *defined by*

$$\gamma := \sup\left\{ \frac{\|h(\phi(t,0,u))\|_Y}{\sup_{0\leq s\leq t}(\|u(s)\|_U)} : t \geq 0, u \in M(U), \sup_{0\leq s\leq t}(\|u(s)\|_U) > 0 \right\} \quad (2.7)$$

*is called the minimum IOS-gain of the input* $u$.

**(d)** *A linear control system* $\Sigma = (X,Y,M(U),\phi,h)$ *is called Exponentially Input-to-State Stable (exp-ISS), if it is exp-IOS with output the identity map* $h(x)=x$. *In this case, the constant* $g \geq 0$ *involved in (2.6) is called an ISS-gain of the input* $u$. *Moreover, if* $U \neq \{0\}$ *then the constant* $\gamma \geq 0$ *defined by (2.7) is called the minimum ISS-gain of the input* $u$.

**(e)** *An OES linear control system* $\Sigma = (X,Y,M(U),\phi,h)$ *satisfies the output asymptotic gain property, if there exists a constant* $g_{as} \geq 0$ *such that for every* $u \in M(U)$ *with* $\sup_{s\geq 0}(\|u(s)\|_U) < +\infty$ *it holds that*

$$\limsup_{t\to+\infty}\left(\|h(\phi(t,0,u))\|_Y\right) \leq g_{as} \sup_{s\geq 0}(\|u(s)\|_U) \quad (2.8)$$

*The constant* $g_{as} \geq 0$ *is called an Output Asymptotic Gain (OAG) of the input* $u$. *If* $U \neq \{0\}$ *then the constant* $\gamma_{as} \geq 0$ *defined by*

$$\gamma_{as} := \sup\left\{ \frac{\limsup_{t\to+\infty}\left(\|h(\phi(t,0,u))\|_Y\right)}{\sup_{s\geq 0}(\|u(s)\|_U)} : u \in M(U), 0 < \sup_{s\geq 0}(\|u(s)\|_U) < +\infty \right\} \quad (2.9)$$

*is called the minimum OAG of the input* $u$.

**(f)** *An ES linear control system* $\Sigma = (X,Y,M(U),\phi,h)$ *satisfies the asymptotic gain property, if it satisfies the output asymptotic gain property with output the identity map* $h(x)=x$. *In this case, the constant* $g_{as} \geq 0$ *involved in (2.8) is called an Asymptotic Gain (AG) of the input* $u$. *Moreover, if* $U \neq \{0\}$ *then the constant* $\gamma_{as} \geq 0$ *defined by (2.9) is called the minimum AG of the input* $u$.

**Remark 2.4:**
   **(i)** The minimum IOS-gain of the input $u$ is an IOS-gain of the input $u$, i.e., $\|h(\phi(t,x,u))\|_Y \leq M\exp(-\sigma t)\|x\|_X + \gamma \sup_{0\leq s\leq t}(\|u(s)\|_U)$ for all $t \geq 0$, $x \in X$, $u \in M(U)$. Similarly, the minimum OAG of the input $u$ is an OAG of the input $u$, i.e., $\limsup_{t\to+\infty}(\|h(\phi(t,0,u))\|_Y) \leq \gamma_{as}\sup_{s\geq 0}(\|u(s)\|_U)$ for every $u \in M(U)$. Furthermore, the following equations hold:

$$\gamma = \inf\{g: g \text{ is an IOS gain for } \Sigma\} \quad (2.10)$$

$$\gamma_{as} = \inf\{g_{as}: g_{as} \text{ is an OAG for } \Sigma\} \quad (2.11)$$

   **(ii)** If $U \neq \{0\}$ then $\gamma_{as} \leq \gamma$.
   **(iii)** Due to linearity the following equations hold:

$$\gamma = \sup\left\{ \|h(\phi(t,0,u))\|_Y : t \geq 0, u \in M(U), \sup_{0 \leq s \leq t}(\|u(s)\|_U) = 1 \right\} \quad (2.12)$$

$$\gamma_{as} = \sup\left\{ \limsup_{t \to +\infty}(\|h(\phi(t,0,u))\|_Y) : u \in M(U), \sup_{s \geq 0}(\|u(s)\|_U) = 1 \right\} \quad (2.13)$$

**(iv)** The inequality $\|h(\phi(t,0,u))\|_Y \leq \gamma \sup_{s \geq 0}(\|u(s)\|_U)$ for all $t \geq 0$, $u \in M(U)$, that follows from the fact that the minimum IOS-gain of the input $u$ is an IOS-gain of the input $u$, implies that $\sup_{t \geq 0}\left\{ \sup(\|h(\phi(t,0,u))\|_Y) : u \in M(U), \sup_{s \geq 0}(\|u(s)\|_U) = 1 \right\} \leq \gamma$. Moreover, since for every $t > 0$, $u \in M(U)$ with $\sup_{0 \leq s \leq t}(\|u(s)\|_U) = 1$ there exists an input $\tilde{u} \in M(U)$ with $\sup_{s \geq 0}(\|\tilde{u}(s)\|_U) = 1$ and $u(s) = \tilde{u}(s)$ for $s \in [0,t]$ (namely, the input $\tilde{u}(s) := u(s)$ for $s \in [0,t]$, $\tilde{u}(s) := 0$ for $s > t$), it follows from causality that for every $t > 0$, $u \in M(U)$ with $\sup_{0 \leq s \leq t}(\|u(s)\|_U) = 1$ the inequality $\|h(\phi(t,0,u))\|_Y \leq \sup_{s \geq 0}\left\{ \sup(\|h(\phi(s,0,v))\|_Y) : v \in M(U), \sup_{s \geq 0}(\|v(s)\|_U) = 1 \right\}$. Using (2.12), the fact that $\|h(\phi(0,0,u))\|_Y = 0$ (a consequence of the identity property and the linearity of the output map) and the previous inequality, we obtain that $\gamma \leq \sup\left\{ \sup_{s \geq 0}(\|h(\phi(s,0,v))\|_Y) : v \in M(U), \sup_{s \geq 0}(\|v(s)\|_U) = 1 \right\}$. Combining, we get:

$$\gamma = \sup\left\{ \sup_{t \geq 0}(\|h(\phi(t,0,u))\|_Y) : u \in M(U), \sup_{s \geq 0}(\|u(s)\|_U) = 1 \right\} \quad (2.14)$$

**(v)** An exp-ISS linear system is also an exp-IOS system for any continuous and linear output map $h: X \to Y$. Indeed, if there exist constants $M, \sigma > 0$, $g \geq 0$ such that for every $(x,u) \in X \times M(U)$ and $t \geq 0$ it holds that $\|\phi(t,x,u)\|_X \leq M \exp(-\sigma t)\|x\|_X + g \sup_{0 \leq s \leq t}(\|u(s)\|_U)$, then the inequality $\|h(\phi(t,x,u))\|_Y \leq KM \exp(-\sigma t)\|x\|_X + Kg \sup_{0 \leq s \leq t}(\|u(s)\|_U)$ will also hold for every $(x,u) \in X \times M(U)$ and $t \geq 0$, where $K := \sup\{\|h(x)\|_Y : x \in X, \|x\|_X = 1\}$ is the norm of $h$. Moreover, if $U \neq \{0\}$ then the minimum IOS-gain is less or equal to $K\gamma$, where $\gamma$ is the minimum ISS-gain.

Finally, we end this section with an additional notion.

**Definition 2.5:** *A linear control system* $\Sigma = (X, Y, M(U), \phi, h)$ *is called robust with respect to input $u$ if $U \neq \{0\}$ and there exists a non-decreasing function $b: \Re_+ \to \Re_+$ such that for every $u \in M(U)$ and $t \geq 0$ it holds that*

$$\|\phi(t,0,u)\|_X \leq b(t) \sup_{0 \leq s \leq t}(\|u(s)\|_U) \quad (2.15)$$

The notion of robustness with respect to an input is related to the notion of admissibility (see [12,21] and the references therein). Clearly, every exp-ISS linear system with $U \neq \{0\}$ is a system which is robust with respect to input $u$ (it satisfies (2.15) with $b(t) \equiv g$, where $g \geq 0$ is any ISS-gain; a direct consequence of (2.6)).

## 3. Motivating Example

This section is devoted to the presentation of a motivating example where our estimates indicate different values for the IOS gain and the OAG. The example is dealing with a linear conventional time-delay system with discrete and distributed delays. The system satisfies the exp-ISS property.

Consider the linear distributed delay system

$$\dot{y}(t) = Ay(t) + Bz(t-\tau) + Gu(t)$$
$$\dot{z}(t) = (KB - \mu I_m)z(t) + K(A + \mu I_n)\left(\exp(A\tau)y(t) + \int_{t-\tau}^{t} \exp(A(t-s))Bz(s)ds\right) \quad (3.1)$$

where $y(t) \in \Re^n$, $z(t) \in \Re^m$, $u(t) \in \Re^p$, $A \in \Re^{n \times n}$, $B \in \Re^{n \times m}$, $G \in \Re^{n \times p}$, $K \in \Re^{m \times n}$ are real matrices for which the matrix $(A+BK) \in \Re^{n \times n}$ is a Hurwitz matrix and $\tau, \mu > 0$ are positive constants. System (3.1) arises in the application of predictor feedback for the linear system with input delay

$$\dot{y}(t) = Ay(t) + Bz(t-\tau) + Gu(t)$$

where $y(t) \in \Re^n$ is the state, $z(t) \in \Re^m$ is the control input and $u(t) \in \Re^p$ is a disturbance. The predictor feedback law $z(t) = K\exp(A\tau)y(t) + K\int_{t-\tau}^{t}\exp(A(t-s))Bz(s)ds$ is implemented dynamically (see pages 23-24 in [6]; dynamic implementation of predictor feedback). The state space of system (1) is $X = \Re^n \times C^0([-\tau,0];\Re^m)$ with state $(y(t), z_t) \in \Re^n \times C^0([-\tau,0];\Re^m)$, where $z_t$ denotes the $\tau$-history of $z$ at time $t \geq 0$, i.e., $(z_t)(\theta) = z(t+\theta)$ for $\theta \in [-\tau, 0]$. For system (3.1) we also have In this case we also have $Y = \Re^n$ with output map $h(y(t), z_t) = y(t)$, while $U = \Re^p$ and $M(U) = L_{loc}^{\infty}(\Re_+; \Re^p)$ is the set of all Lebesgue measurable and locally essentially bounded inputs $u : \Re_+ \to \Re^p$. By virtue of linearity the solution of system (3.1) with arbitrary initial condition $(y(0), z_0) \in \Re^n \times C^0([-\tau,0];\Re^m)$ exists for all $t \geq 0$.

We are interested in studying the robustness properties of system (3.1) and more specifically the effect of the input $u$ to the state component $y$. To this purpose, we notice that the following differential equation holds for $\xi(t) = z(t) - K\exp(A\tau)y(t) - K\int_{t-\tau}^{t}\exp(A(t-s))Bz(s)ds$:

$$\dot{\xi}(t) = -\mu\xi(t) - K\exp(A\tau)Gu(t), \text{ for } t \geq 0 \text{ a.e.} \quad (3.2)$$

Using the variations of constants formula and solving (3.2), we obtain the following formula for $t \geq 0$:

$$\xi(t) = \exp(-\mu t)\xi(0) - \int_{0}^{t}\exp(-\mu(t-s))K\exp(A\tau)Gu(s)ds \quad (3.3)$$

Combining (3.1) and (3.3), definition $\xi(t) = z(t) - K\exp(A\tau)y(t) - K\int_{t-\tau}^{t}\exp(A(t-s))Bz(s)ds$ and the fact that $y(t+\tau) = \exp(A\tau)y(t) + \int_{t-\tau}^{t}\exp(A(t-s))Bz(s)ds + \int_{t}^{t+\tau}\exp(A(t+\tau-s))Gu(s)ds$, we get the differential equation for $t \geq 0$ a.e.:

$$\begin{aligned}\dot{y}(t+\tau) &= (A+BK)y(t+\tau) + B\exp(-\mu t)\xi(0) \\ &+ Gu(t+\tau) - \int_{t}^{t+\tau} BK\exp(A(t+\tau-s))Gu(s)ds - \int_{0}^{t}\exp(-\mu(t-s))BK\exp(A\tau)Gu(s)ds\end{aligned} \quad (3.4)$$

Using the variations of constants formula and solving (3.1) and (3.4), we obtain the following formulas

$$y(t) = \exp(At)y(0) + \int_{0}^{t}\exp(A(t-s))Bz(s-\tau)ds + P(t)u, \text{ for } t \in [0,\tau] \quad (3.5)$$

$$\begin{aligned}y(t) &= Q(t)\left(\exp(A\tau)y(0) + \int_{-\tau}^{0}\exp(-As)Bz(s)ds\right) \\ &+ \int_{0}^{t-\tau}\exp((A+BK)(t-\tau-s))\exp(-\mu s)dsBz(0) + P(t)u\end{aligned}$$

$$\text{for } t \geq \tau \quad (3.6)$$

where $Q(t) := \exp((A+BK)(t-\tau))\left(I - \int_{0}^{t-\tau}\exp(-(A+BK)s)\exp(-\mu s)dsBK\right)$ and $P(t)$ is the operator defined by:

$$P(t)u = \int_{0}^{t}\exp(A(t-s))Gu(s)ds, \text{ for } t \in [0,\tau] \quad (3.7)$$

$$\begin{aligned}P(t)u &= \exp((A+BK)(t-\tau))\int_{0}^{\tau}\exp(A(\tau-s))Gu(s)ds \\ &+ \int_{\tau}^{t}\exp((A+BK)(t-s))\left(Gu(s) - \int_{s-\tau}^{s}BK\exp(A(s-l))Gu(l)dl\right)ds \\ &- \int_{0}^{t-\tau}\exp((A+BK)(t-\tau-s))\left(\int_{0}^{s}\exp(-\mu(s-l))BK\exp(A\tau)Gu(l)dl\right)ds\end{aligned}$$

$$\text{for } t \geq \tau \quad (3.8)$$

Since $(A+BK) \in \mathfrak{R}^{n \times n}$ is a Hurwitz matrix, by using (3.3), (3.5), (3.6), (3.7), (3.8) and the fact that $z(t) = Ky(t+\tau) + \xi(t) - K\int_{t}^{t+\tau}\exp(A(t+\tau-s))Gu(s)ds$ for $t \geq 0$, we obtain an ISS estimate of the form

$$|y(t)| + \|z_t\| \leq G\exp(-\sigma t)(|y(0)| + \|z_0\|) + K\sup_{0 \leq s \leq t}(|u(s)|), \text{ for } t \geq 0 \quad (3.9)$$

where $G, \sigma, K > 0$ are constants (independent of $t, y(0), z_0$ and the input $u$) and $\|z_t\| = \max_{\theta \in [-\tau, 0]} (|z(t+\theta)|)$. Moreover, using (3.5), (3.6), (3.7), (3.8), the fact that there exist constants $M \geq 1$, $\sigma > 0$ such that $|\exp((A+BK)t)| \leq M \exp(-\sigma t)$ for $t \geq 0$, we obtain the following estimate for $t \geq 0$

$$|y(t)| \leq R \exp(-\sigma t)(|y(0)| + \|z_0\|)$$
$$+ M \sigma^{-1} \left( |G| + \int_0^\tau \phi(s) dl + \mu^{-1} \phi(\tau) + \sigma \min(1, \exp(-\sigma(t-\tau))) \int_0^\tau r(s) ds \right) \sup_{0 \leq s \leq t} (|u(s)|) \quad (3.10)$$

where $R > 0$ is an appropriate constant and $\phi, r : \mathfrak{R}_+ \to \mathfrak{R}_+$ are continuous functions that satisfy $|BK \exp(As)G| \leq \phi(s)$ and $|\exp(As)G| \leq r(s)$ for $s \geq 0$. It follows from (3.10) that an IOS-gain of the input $u$ is $M \sigma^{-1} \left( |G| + \int_0^\tau \phi(s) dl + \mu^{-1} \phi(\tau) + \sigma \int_0^\tau r(s) ds \right)$. On the other hand, estimate (3.10) implies the following inequality for every bounded input $u$

$$\limsup_{t \to +\infty} (|y(t)|) \leq M \sigma^{-1} \left( |G| + \int_0^\tau \phi(s) dl + \mu^{-1} \phi(\tau) \right) \sup_{s \geq 0} (|u(s)|) \quad (3.11)$$

which shows that an OAG of the input $u$ is $M \sigma^{-1} \left( |G| + \int_0^\tau \phi(s) dl + \mu^{-1} \phi(\tau) \right)$. Therefore, in this case our estimates indicate that the minimum OAG is less than the minimum IOS-gain of the input $u$.

The results of the following section show that this is not the case. For a linear system like (3.1), the minimum OAG is equal to the minimum IOS-gain. Therefore, for this example, we are in a position to guarantee that the minimum IOS-gain is less or equal to $M \sigma^{-1} \left( |G| + \int_0^\tau \phi(s) dl + \mu^{-1} \phi(\tau) \right)$.

## 4. Main Results

Our first main result guarantees the exp-ISS property and provides a useful upper bound for the minimum ISS gain. In order to be able to state our first main result, we notice that for an ES linear system $\Sigma = (X, Y, M(U), \phi, h)$, for each $t \geq 0$ there exists a linear continuous operator $S(t) : X \to X$ such that $\phi(t, x, 0) = S(t)x$. The one-parameter family of linear bounded operators is a semigroup of bounded linear operators on $X$ (see Definition 1.1 in [13]), which (by virtue of (2.5)) satisfies the estimate $\|S(t)\|_{L(X)} \leq M \exp(-\sigma t)$ for all $t \geq 0$, where $L(X)$ is the space of bounded linear operators $A : X \to X$. It should be noticed here that *we do not assume anything* for the continuity with respect to $t \geq 0$ of the semigroup $S(t) : X \to X$.

**Theorem 4.1:** *If a linear control system $\Sigma = (X, Y, M(U), \phi, h)$ is robust with respect to input $u$ and ES then it is exp-ISS. Moreover, the following formula holds for the minimum ISS-gain*

$$\gamma \leq \min \left( \inf \left\{ \sup_{0 \leq s < T} \left( \frac{M(\sigma) \exp(-\sigma s) b(T)}{1 - M(\sigma) \exp(-\sigma T)} + b(s) \right) : \sigma \in F_S, T > \frac{1}{\sigma} \ln(M(\sigma)) \right\}, \sup_{s \geq 0} (b(s)) \right) \quad (4.1)$$

where

$$M(\sigma) := \sup_{t \geq 0} \left( \|S(t)\|_{L(X)} \exp(\sigma t) \right) \quad (4.2)$$

$$F_S := \{\sigma > 0 : M(\sigma) < +\infty\} \quad (4.3)$$

and $b : \Re_+ \to \Re_+$ is the non-decreasing function involved in (2.15).

**Remark 4.2:** The converse of Theorem 4.1 also holds, i.e., an exp-ISS linear system is an ES system which is robust with respect to input $u$. Notice that (4.1) in conjunction with definitions (4.2) and (4.3) imply that $\gamma \leq \sup_{0 \leq s < T} \left( \frac{M \exp(-\sigma s) b(T)}{1 - M \exp(-\sigma T)} + b(s) \right) \leq \left( \frac{1 + M(1 - \exp(-\sigma T))}{1 - M \exp(-\sigma T)} \right) b(T) < +\infty$, where $M, \sigma > 0$ are any constants for which estimate (2.5) holds and $T > \frac{1}{\sigma} \ln(M)$.

In order to state the next main result, we need to define for a linear control system $\Sigma = (X, Y, M(U), \phi, h)$ with $U \neq \{0\}$ the class of periodic inputs, namely,

$$M_{per}(U) := \{ u \in M(U) : u \text{ periodic} \} \quad (4.4)$$

Moreover, for every $T > 0$, we introduce the following optimal control problem:

$$V(T) := \sup \left\{ \|h(\phi(T, 0, u))\|_Y : u \in M(U), \sup_{0 \leq t \leq T} \left( \|u(t)\|_U \right) = 1 \right\} \quad (4.5)$$

It is clear that if $\Sigma = (X, Y, M(U), \phi, h)$ is robust with respect to input $u$ then the optimal control problem (4.5) is well-defined for every $T > 0$ and the function $V : (0, +\infty) \to \Re_+$ has finite values. We are now in a position to state our second main result.

**Theorem 4.3:** *Suppose that a linear control system $\Sigma = (X, Y, M(U), \phi, h)$ with $U \neq \{0\}$ is robust with respect to input $u$ and ES. Moreover, suppose that the following property holds:*

**Strict Causality:** *for every $t > 0$ and for every $u, \tilde{u} \in M(U)$ with $\tilde{u}(s) = u(s)$ for all $s \in [0, t)$, it holds that $\phi(t, 0, u) = \phi(t, 0, \tilde{u})$.*

*Then the following equations hold:*

$$\gamma = \gamma_{as} = \sup_{T > 0} (V(T)) \quad (4.6)$$

$$\gamma := \sup \left\{ \sup_{t \geq 0} \left( \|h(\phi(t, 0, u))\|_Y \right) : u \in M_{per}(U), \sup_{s \geq 0} \left( \|u(s)\|_U \right) = 1 \right\} \quad (4.7)$$

$$\gamma_{as} := \sup \left\{ \limsup_{t \to +\infty} \left( \|h(\phi(t, 0, u))\|_Y \right) : u \in M_{per}(U), \sup_{s \geq 0} \left( \|u(s)\|_U \right) = 1 \right\} \quad (4.8)$$

The reader should notice the difference between strict causality and causality in (iv)(c) of Definition 2.1: in property (iv)(c) of Definition 2.1, the inputs $u, \tilde{u} \in M(U)$ for which $\phi(t, 0, u) = \phi(t, 0, \tilde{u})$ must also satisfy $u(t) = \tilde{u}(t)$.

Theorem 4.3 is important because in cases like the example of the previous section, it can allow a more accurate estimation of the IOS-gain (by estimating the OAG). The fact that the minimum IOS-gain and the OAG can be estimated using only periodic inputs is important: in the finite-dimensional case a periodic input creates a periodic solution which attracts exponentially all solutions. Formulas (4.6), (4.7), (4.8) can allow an explicit and exact calculation of the minimum IOS-gain (see Section 5). Finally, in many cases formulas (4.6), (4.7), (4.8) can allow the derivation of upper and lower bounds for the minimum IOS-gain.

We next proceed to the proofs of the two main results.

**Proof of Theorem 4.1:** The fact that $\gamma \leq \sup_{s \geq 0}(b(s))$ is a direct consequence of (2.14) and inequality (2.15). Therefore, we only have to prove that $\gamma \leq \inf\left\{\sup_{0 \leq s < T}\left(\frac{M(\sigma)\exp(-\sigma s)b(T)}{1-M(\sigma)\exp(-\sigma T)} + b(s)\right) : \sigma \in F_S, T > \frac{1}{\sigma}\ln(M(\sigma))\right\}$. To this purpose, it suffices to show that for every pair of constants $M, \sigma > 0$ for which estimate (2.5) holds, the following estimate holds:

$$\gamma \leq \sup_{0 \leq s < T}\left(\frac{M\exp(-\sigma s)b(T)}{1-M\exp(-\sigma T)} + b(s)\right), \text{ for all } T > \frac{1}{\sigma}\ln(M) \tag{4.9}$$

Let $M, \sigma > 0$ be a pair of constants for which estimate (2.5) holds and let $T > \frac{1}{\sigma}\ln(M)$ be given (arbitrary). Define

$$\lambda := M\exp(-\sigma T) < 1 \tag{4.10}$$

Using (2.5), (2.15), (4.10), the semigroup property and the triangle inequality, we obtain for all integers $k \geq 0$ and all $(x, u) \in X \times M(U)$:

$$\begin{aligned}
\|\phi((k+1)T, x, u)\|_X &= \|\phi(T, \phi(kT, x, u), \delta_{kT}u)\|_X \\
&= \|\phi(T, \phi(kT, x, u), 0) + \phi(T, 0, \delta_{kT}u)\|_X \\
&\leq \|\phi(T, \phi(kT, x, u), 0)\|_X + \|\phi(T, 0, \delta_{kT}u)\|_X \\
&\leq \lambda\|\phi(kT, x, u)\|_X + b(T)\sup_{0 \leq s \leq T}\left(\|(\delta_{kT}u)(s)\|_U\right) \\
&\leq \lambda\|\phi(kT, x, u)\|_X + b(T)\sup_{kT \leq s \leq (k+1)T}\left(\|u(s)\|_U\right)
\end{aligned} \tag{4.11}$$

Using induction and (4.11), we guarantee that the following estimate holds for all integers $k \geq 1$ and all $(x, u) \in X \times M(U)$:

$$\|\phi(kT, x, u)\|_X \leq \lambda^k \|x\|_X + b(T)\sum_{j=0}^{k-1}\lambda^{k-1-j}\sup_{jT \leq s \leq (j+1)T}\left(\|u(s)\|_U\right) \tag{4.12}$$

By virtue of (4.12) and the fact that $\lambda \in (0,1)$, the following estimates hold for all integers $k \geq 1$ and all $(x, u) \in X \times M(U)$:

$$\|\phi(kT,x,u)\|_X \leq \lambda^k \|x\|_X + b(T)\left(\sum_{j=0}^{k-1} \lambda^{k-1-j}\right) \sup_{0 \leq s \leq kT}\left(\|u(s)\|_U\right)$$

$$= \lambda^k \|x\|_X + b(T)\frac{1-\lambda^k}{1-\lambda} \sup_{0 \leq s \leq kT}\left(\|u(s)\|_U\right) \quad (4.13)$$

$$\leq \lambda^k \|x\|_X + \frac{b(T)}{1-\lambda} \sup_{0 \leq s \leq kT}\left(\|u(s)\|_U\right)$$

The identity property guarantees that (4.13) holds also for $k=0$, i.e., the following estimate holds for all integers $k \geq 0$ and all $(x,u) \in X \times M(U)$:

$$\|\phi(kT,x,u)\|_X \leq \lambda^k \|x\|_X + \frac{b(T)}{1-\lambda} \sup_{0 \leq s \leq kT}\left(\|u(s)\|_U\right) \quad (4.14)$$

Let arbitrary $t \geq 0$ and let $k$ be the integer part of $t/T$. Using (2.5), (2.15), (4.14), the triangle inequality and the semigroup property, we get for all $(x,u) \in X \times M(U)$:

$$\begin{aligned}
\|\phi(t,x,u)\|_X &= \|\phi(t-kT,\phi(kT,x,u),\delta_{kT}u)\|_X \\
&= \|\phi(t-kT,\phi(kT,x,u),0) + \phi(t-kT,0,\delta_{kT}u)\|_X \\
&\leq \|\phi(t-kT,\phi(kT,x,u),0)\|_X + \|\phi(t-kT,0,\delta_{kT}u)\|_X \\
&\leq M\exp(-\sigma(t-kT))\|\phi(kT,x,u)\|_X + b(t-kT)\sup_{0 \leq s \leq t-kT}\left(\|(\delta_{kT}u)(s)\|_U\right) \quad (4.15) \\
&\leq M\exp(-\sigma(t-kT))\|\phi(kT,x,u)\|_X + b(t-kT)\sup_{kT \leq s \leq t}\left(\|u(s)\|_U\right) \\
&\leq M\exp(-\sigma(t-kT))\lambda^k \|x\|_X + \left(\exp(-\sigma(t-kT))\frac{Mb(T)}{1-\lambda} + b(t-kT)\right)\sup_{0 \leq s \leq t}\left(\|u(s)\|_U\right)
\end{aligned}$$

Using the fact that $0 \leq t-kT < T$, definition (4.10) and (4.15), we get for all $(x,u) \in X \times M(U)$:

$$\begin{aligned}
\|\phi(t,x,u)\|_X &\leq M^{k+1}\exp(-\sigma t)\|x\|_X + \left(\frac{M\exp(-\sigma(t-kT))b(T)}{1-M\exp(-\sigma T)} + b(t-kT)\right)\sup_{0 \leq s \leq t}\left(\|u(s)\|_U\right) \\
&= M\exp\left(-\sigma t + k\ln(M)\right)\|x\|_X + \sup_{0 \leq r < T}\left(\frac{M\exp(-\sigma r)b(T)}{1-M\exp(-\sigma T)} + b(r)\right)\sup_{0 \leq s \leq t}\left(\|u(s)\|_U\right)
\end{aligned} \quad (4.16)$$

Using the fact that $k \leq t/T$ and $M \geq 1$ (a direct consequence of (2.5) and the identity property), we obtain from (4.16) for all $(x,u) \in X \times M(U)$:

$$\|\phi(t,x,u)\|_X \leq M\exp\left(-\left(\sigma - T^{-1}\ln(M)\right)t\right)\|x\|_X + \sup_{0 \leq r < T}\left(\frac{M\exp(-\sigma r)b(T)}{1-M\exp(-\sigma T)} + b(r)\right)\sup_{0 \leq s \leq t}\left(\|u(s)\|_U\right) \quad (4.17)$$

Since $T > \frac{1}{\sigma}\ln(M)$, it follows that estimate (4.17) shows the exp-ISS property with ISS-gain $\sup_{0 \leq r < T}\left(\frac{M\exp(-\sigma r)b(T)}{1-M\exp(-\sigma T)} + b(r)\right)$. Therefore, (4.9) holds. The proof is complete. ◁

**Proof of Theorem 4.3:** Definitions (2.12), (4.5) and the fact that $\gamma_{as} \leq \gamma$ imply that the following inequalities hold:

$$\gamma_{as} \leq \gamma \leq \sup_{T>0}(V(T)) \qquad (4.18)$$

Consequently, in order to prove (4.6) it suffices to show that $\gamma_{as} \geq \sup_{T>0}(V(T))$. Taking into account (2.13), it suffices to show that for every $\varepsilon > 0$ there exists an input $u_\varepsilon \in M(U)$ with $\sup_{s \geq 0}(\|u_\varepsilon(s)\|_U) = 1$ such that $\limsup_{t \to +\infty}(\|h(\phi(t,0,u_\varepsilon))\|_Y) \geq \sup_{T>0}(V(T)) - \varepsilon$.

Let $\varepsilon > 0$ be given (arbitrary). There exists $T := T(\varepsilon) > 0$ such that

$$V(T) \geq \sup_{t>0}(V(t)) - \varepsilon/3 \qquad (4.19)$$

Moreover, it follows from definition (4.5) that there exists $w_\varepsilon \in M(U)$ with $\sup_{0 \leq t \leq T}(\|w_\varepsilon(t)\|_U) = 1$ such that

$$\|h(\phi(T,0,w_\varepsilon))\|_Y \geq V(T) - \varepsilon/3 \qquad (4.20)$$

Since $\Sigma = (X, Y, M(U), \phi, h)$ is ES, exist constants $M, \sigma > 0$ such that for every $x \in X$ and $t \geq 0$ (2.5) holds. By continuity of $h: X \to Y$, there exists $\delta := \delta(\varepsilon) > 0$ such that for all $\xi \in X$ with $\|\xi - \phi(T,0,w_\varepsilon)\|_X \leq \delta$ the following inequality holds:

$$\|h(\xi) - h(\phi(T,0,w_\varepsilon))\|_Y \leq \varepsilon/3 \qquad (4.21)$$

Since $\Sigma = (X, Y, M(U), \phi, h)$ with $U \neq \{0\}$ is robust with respect to input $u$, there exists a non-decreasing function $b: \Re_+ \to \Re_+$ such that for every $u \in M(U)$ and $t \geq 0$, inequality (2.15) holds.

Let $R := R(\varepsilon) > 0$ be sufficiently large so that

$$M \exp(-\sigma R)\left(\delta \exp(-\sigma) + b(T+1)\right) \leq M^{-1} \delta \exp(\sigma T) \qquad (4.22)$$

and define

$$v_\varepsilon(t) = \begin{cases} w_\varepsilon(t), & t \in [0, T] \\ 0, & t > T \end{cases} \qquad (4.23)$$

Since $0 \in M(U)$ (notice that here denotes the zero input), property (iv) of Definition 2.1 guarantees $v_\varepsilon$ as defined by (4.23) is an input of class $M(U)$. Moreover, notice that property (v) of Definition (2.2) guarantees that its $(1+T+R)$-periodic extension of $v_\varepsilon$ given by

$$u_\varepsilon(t) = v_\varepsilon\left(t - (T+R+1)\left\lfloor\frac{t}{T+R+1}\right\rfloor\right), \text{ for } t \geq 0 \qquad (4.24)$$

is a function of class $M(U)$. Definitions (4.23), (4.24) imply that $\sup_{t \geq 0}(\|u_\varepsilon(t)\|_U) = \sup_{t \geq 0}(\|v_\varepsilon(t)\|_U) = \sup_{0 \leq t \leq T}(\|w_\varepsilon(t)\|_U) = 1$. Moreover, it follows from (2.15) and the fact

that $\sup_{t\geq 0}\left(\left\|\delta_\tau u_\varepsilon(t)\right\|_U\right)=1$ for all $\tau>0$ (recall that $\delta_\tau$ for $\tau>0$ is the shift operator that satisfies $(\delta_\tau u)(t):=u(t+\tau)$ for $t\geq 0$ and every $u\in M(U)$) that

$$\left\|\phi(T+1,0,\delta_\tau u_\varepsilon)\right\|_X \leq b(T+1), \text{ for all } \tau>0 \tag{4.25}$$

We next show by induction on $k\geq 0$ that

$$\left\|\phi(k(T+R),0,u_\varepsilon)\right\|_X \leq M^{-1}\delta\exp(\sigma T), \text{ for all } k\geq 0 \tag{4.26}$$

Indeed, the identity property guarantees (4.26) for $k=0$. Assuming that (4.26) holds for certain $k\geq 0$, we get from the semigroup property:

$$\begin{aligned}&\left\|\phi((k+1)(T+R+1),0,u_\varepsilon)\right\|_X = \\ &\left\|\phi(R,\phi((k+1)(T+R+1)-R,0,u_\varepsilon),\delta_{(k+1)(T+R+1)-R}u_\varepsilon)\right\|_X\end{aligned} \tag{4.27}$$

The input $\delta_{(k+1)(T+R+1)-R}u_\varepsilon$ satisfies $(\delta_{(k+1)(T+R+1)-R}u_\varepsilon)(s)=0$ for $s\in[0,R)$ (recall definitions (4.23), (4.24)). It follows from the strict causality property and (2.5), (4.27) that

$$\begin{aligned}&\left\|\phi((k+1)(T+R+1),0,u_\varepsilon)\right\|_X = \\ &\left\|\phi(R,\phi((k+1)(T+R+1)-R,0,u_\varepsilon),0)\right\|_X \\ &\leq M\exp(-\sigma R)\left\|\phi((k+1)(T+R+1)-R,0,u_\varepsilon)\right\|_X\end{aligned} \tag{4.28}$$

Using again the semigroup property, the linearity of the response map and the triangle inequality, we get from (4.28):

$$\begin{aligned}&\left\|\phi((k+1)(T+R+1),0,u_\varepsilon)\right\|_X \\ &\leq M\exp(-\sigma R)\left\|\phi(T+1,\phi(k(T+R+1),0,u_\varepsilon),\delta_{k(T+R+1)}u_\varepsilon)\right\|_X \\ &\leq M\exp(-\sigma R)\left\|\phi(T+1,\phi(k(T+R+1),0,u_\varepsilon),0)\right\|_X \\ &+M\exp(-\sigma R)\left\|\phi(T+1,0,\delta_{k(T+R+1)}u_\varepsilon)\right\|_X\end{aligned} \tag{4.29}$$

Using (4.29) in conjunction with (4.25) and (2.5), we get:

$$\begin{aligned}&\left\|\phi((k+1)(T+R+1),0,u_\varepsilon)\right\|_X \\ &\leq M^2\exp(-\sigma(T+R+1))\left\|\phi(k(T+R+1),0,u_\varepsilon)\right\|_X + M\exp(-\sigma R)b(T+1)\end{aligned} \tag{4.30}$$

Estimate (4.30) in conjunction with (4.26) gives:

$$\left\|\phi((k+1)(T+R+1),0,u_\varepsilon)\right\|_X \leq M\exp(-\sigma R)\left(\delta\exp(-\sigma)+b(T+1)\right) \tag{4.31}$$

Inequality (4.26) with $k+1$ in place of $k$ is a direct consequence of (4.31) and (4.22).

We next evaluate the quantity $\left\|\phi(T+k(T+R+1),0,u_\varepsilon)-\phi(T,0,w_\varepsilon)\right\|_X$ for all $k\geq 0$. Using the semigroup property and the fact that $(\delta_{k(T+R+1)}u_\varepsilon)(s)=w_\varepsilon(s)$ for all $s\in[0,T]$, the causality property and the linearity of the response map, we get for all $k\geq 0$:

$$\begin{aligned}
&\|\phi(T+k(T+R+1),0,u_\varepsilon) - \phi(T,0,w_\varepsilon)\|_X \\
&= \|\phi(T,\phi(k(T+R+1),0,u_\varepsilon),\delta_{k(T+R+1)}u_\varepsilon) - \phi(T,0,w_\varepsilon)\|_X \\
&= \|\phi(T,\phi(k(T+R+1),0,u_\varepsilon),w_\varepsilon) - \phi(T,0,w_\varepsilon)\|_X \\
&= \|\phi(T,\phi(k(T+R+1),0,u_\varepsilon),0)\|_X
\end{aligned} \qquad (4.32)$$

It follows from (4.32), (2.5) and (4.26) that the following estimate holds for all $k \geq 0$:

$$\begin{aligned}
&\|\phi(T+k(T+R+1),0,u_\varepsilon) - \phi(T,0,w_\varepsilon)\|_X \\
&\leq M \exp(-\sigma T) \|\phi(k(T+R+1),0,u_\varepsilon)\|_X \leq \delta
\end{aligned} \qquad (4.33)$$

Consequently, (4.21) in conjunction with (4.33) give for all $k \geq 0$:

$$\|h(\phi(T+k(T+R+1),0,u_\varepsilon)) - h(\phi(T,0,w_\varepsilon))\|_Y \leq \varepsilon/3 \qquad (4.34)$$

Therefore, we get from (4.34), (4.19), (4.20) for all $k \geq 0$:

$$\begin{aligned}
\|h(\phi(T+k(T+R+1),0,u_\varepsilon))\|_Y &\geq \|h(\phi(T,0,w_\varepsilon))\|_Y - \varepsilon/3 \\
&\geq V(T) - 2\varepsilon/3 \geq \sup_{t>0}(V(t)) - \varepsilon
\end{aligned} \qquad (4.35)$$

Inequality (4.35) implies that $\limsup_{t\to+\infty}\left(\|h(\phi(t,0,u_\varepsilon))\|_Y\right) \geq \sup_{t>0}(V(t)) - \varepsilon$.

It should be noticed that the input $u_\varepsilon$ constructed above, is a function of class $M_{per}(U)$ with $\sup_{s\geq 0}\left(\|u_\varepsilon(s)\|_U\right) = 1$. Since (4.6) holds, it follows that for every $\varepsilon > 0$ there exists an input $u_\varepsilon \in M_{per}(U)$ with $\sup_{s\geq 0}\left(\|u_\varepsilon(s)\|_U\right) = 1$ such that $\limsup_{t\to+\infty}\left(\|h(\phi(t,0,u_\varepsilon))\|_Y\right) \geq \gamma_{as} - \varepsilon$. Consequently, we get

$$\sup\left\{\limsup_{t\to+\infty}\left(\|h(\phi(t,0,u))\|_Y\right): u \in M_{per}(U), \sup_{s\geq 0}\left(\|u(s)\|_U\right) = 1\right\} \geq \gamma_{as} \qquad (4.36)$$

Equation (2.14) and the fact that $M_{per}(U) \subseteq M(U)$ implies that

$$\sup\left\{\sup_{t\geq 0}\left(\|h(\phi(t,0,u))\|_Y\right): u \in M_{per}(U), \sup_{s\geq 0}\left(\|u(s)\|_U\right) = 1\right\} \leq \gamma \qquad (4.37)$$

Finally, since $\limsup_{t\to+\infty}\left(\|h(\phi(t,0,u))\|_Y\right) \leq \sup_{t\geq 0}\left(\|h(\phi(t,0,u))\|_Y\right)$ for all $u \in M(U)$, the following inequality holds:

$$\begin{aligned}
&\sup\left\{\limsup_{t\to+\infty}\left(\|h(\phi(t,0,u))\|_Y\right): u \in M_{per}(U), \sup_{s\geq 0}\left(\|u(s)\|_U\right) = 1\right\} \\
&\leq \sup\left\{\sup_{t\geq 0}\left(\|h(\phi(t,0,u))\|_Y\right): u \in M_{per}(U), \sup_{s\geq 0}\left(\|u(s)\|_U\right) = 1\right\}
\end{aligned} \qquad (4.38)$$

Inequalities (4.7), (4.8) are direct consequences of equation (4.6) and inequalities (4.36), (4.37), (4.38). The proof is complete. ◁

# 5. Applications to Finite-Dimensional Systems

In this section we focus on finite-dimensional linear time-invariant systems of the form

$$\dot{x} = Ax + Bu, \quad y = Cx, \tag{5.1}$$

where $x \in \mathfrak{R}^n$, $u \in \mathfrak{R}^m$, $y \in \mathfrak{R}^p$ and $A \in \mathfrak{R}^{n \times n}$, $B \in \mathfrak{R}^{n \times m}$, $C \in \mathfrak{R}^{p \times n}$. We assume that $A \in \mathfrak{R}^{n \times n}$ is a Hurwitz matrix so that the assumptions of Theorem 4.1 and Theorem 4.3 automatically hold. It should be noticed that for system (5.1), the optimal control problem (4.5) takes the form for every $T > 0$:

$$\begin{aligned} V(T) &= \sup\left(|Cx(T)|\right) \\ &\quad subject\ to \\ \dot{x}(t) &= Ax(t) + Bu(t),\ t \in [0,T]\ a.e. \\ u &\in L^\infty([0,T]; \mathfrak{R}^m) \\ |u(t)| &\leq 1,\ t \in [0,T]\ a.e. \end{aligned} \tag{5.2}$$

*5.1. Single-Input-Single-Output Systems*

When $m = p = 1$ then system (5.1) is a Single-Input-Single-Output (SISO) system. In this case, the optimal control problem (4.5) becomes for every $T > 0$

$$V(T) = \sup\left\{ \left| \int_0^T C\exp(A(T-s))Bu(s)ds \right| : u \in L^\infty_{loc}(\mathfrak{R}_+; \mathfrak{R}),\ \sup_{0 \leq t \leq T}\left(|u(t)|\right) = 1 \right\} \tag{5.3}$$

Since the function $[0,T] \ni s \to C\exp(A(T-s))B$ is an analytic function, it follows that either is identically zero or vanishes at most a finite number of times. Therefore, the bang-bang control

$$u(s) = \mathrm{sgn}\left(C\exp(A(T-s))B\right),\ s \in [0,T] \tag{5.4}$$

where $\mathrm{sgn}(x) = 1$ for $x \geq 0$ and $\mathrm{sgn}(x) = -1$ for $x < 0$ is a piecewise continuous function on $[0,T]$. Moreover, the bang-bang control given by (5.4) guarantees that

$$\left| \int_0^T C\exp(A(T-s))Bu(s)ds \right| = \int_0^T |C\exp(A(T-s))B|ds \tag{5.5}$$

Since (5.3) implies that $V(T) \leq \int_0^T |C\exp(A(T-s))B|ds$ for all $T > 0$, it follows from (5.5) that

$$V(T) = \int_0^T |C\exp(A(T-s))B|ds \tag{5.6}$$

Consequently, Theorem 4.3 allows us to give an explicit and exact formula for the minimum IOS-gain of system (5.1). More specifically, we obtain the formula:

$$\gamma = \gamma_{as} = \sup_{T>0}(V(T)) = \int_0^{+\infty} |C\exp(As)B|\,ds \tag{5.7}$$

In other words, the minimum IOS-gain is exactly equal to the $L^1$ norm of the output of the input-free system (5.1) (i.e., system (5.1) with $u \equiv 0$) with initial condition $x(0) = B$.

*5.2. A Class of Single-Input Systems*

A special class of single input systems (5.1) for which the optimal control problem (5.2) is explicitly solvable is the class for which the following inequality holds:

$$B'\exp(A't)C'C\exp(As)B \geq 0, \forall t, s \geq 0 \tag{5.8}$$

Indeed, for this class the constant input $u(t) \equiv 1$ solves the optimal control problem (5.2) for every $T > 0$, since we have for every $u \in L^\infty([0,T];\Re)$ with $|u(t)| \leq 1$ for $t \in [0,T]$ a.e.:

$$\begin{aligned}|Cx(T)|^2 &= \int_0^T\int_0^T B'\exp(A'(T-t))C'C\exp(A(T-s))Bu(s)u(t)\,ds\,dt \\ &\leq \int_0^T\int_0^T B'\exp(A'(T-t))C'C\exp(A(T-s))B\,ds\,dt = \left|\int_0^T C\exp(As)B\,ds\right|^2\end{aligned} \tag{5.9}$$

Consequently, equation (4.6) gives:

$$\gamma = \gamma_{as} = \sup_{T>0}\left(\left|\int_0^T C\exp(As)B\,ds\right|\right) \tag{5.10}$$

Since

$$\sup_{T>0}\left(\left|\int_0^T C\exp(As)B\,ds\right|^2\right) = \int_0^{+\infty}\int_0^{+\infty} B'\exp(A't)C'C\exp(As)B\,ds\,dt = \left|\int_0^{+\infty} C\exp(As)B\,ds\right|^2$$

it follows from (5.10) that

$$\gamma = \gamma_{as} = \left|\int_0^{+\infty} C\exp(As)B\,ds\right| = |CA^{-1}B| \tag{5.11}$$

Notice that formula (5.11) does not require the computation of the exponential matrix $\exp(At)$.

As an example, we next give the computation of the minimum ISS-gain for (5.1) under the following assumption:

**(H)** *There exist constants $\lambda_i, q_i > 0$ ($i = 1,...,n$) and an invertible matrix $P \in \Re^{n \times n}$ such that $PAP^{-1} = -diag(\lambda_1,...,\lambda_n)$ and $PP' = diag(q_1,...,q_n)$.*

It is clear that assumption (H) implies that all eigenvalues of $A \in \Re^{n \times n}$ are real. In this case, we have $\exp(At) = P^{-1} diag(\exp(-\lambda_1 t),...,\exp(-\lambda_n t)) P$ and $P^{-1} = P' diag(q_1^{-1},...,q_n^{-1})$. Consequently, (5.8) holds. Using (5.11) with $C = I_n$, we get:

$$\gamma = \gamma_{as} = \sqrt{B'P' diag\left(\frac{1}{q_1 \lambda_1^2},...,\frac{1}{q_n \lambda_n^2}\right) PB} \tag{5.12}$$

As another example of a linear system (5.1) that satisfies inequality (5.8), we can give the case where $A \in \Re^{n \times n}$ is a Metzler matrix and $B \in \Re^{n \times m}$, $C \in \Re^{p \times n}$ are non-negative matrices. In this case, it holds that the vector $C \exp(At) B$ has non-negative entries for every $t \geq 0$ (see Theorem 4.2 on page 64 in [22]). Consequently, inequality (5.8) holds. The case where $A \in \Re^{n \times n}$ is a Metzler matrix and $B \in \Re^{n \times m}$, $C \in \Re^{p \times n}$ are non-negative matrices is important in practice because many physical, chemical and biological systems have non-negative values for their state variables. Therefore, if a linear system (5.1) is used as a model for a system with non-negative states then $A \in \Re^{n \times n}$ must necessarily be a Metzler matrix and $B \in \Re^{n \times m}$, $C \in \Re^{p \times n}$ must be non-negative matrices (see Chapter 4 in [22]).

*5.3. Lower Bounds for the Minimum Gains*

For a single-input system, when a sinusoid input is applied, i.e.,

$$u(t) = \sin(\omega t + \varphi), \text{ for } t \geq 0 \tag{5.13}$$

where $\omega > 0$ and $\varphi \in \Re$, then (for arbitrary initial condition) the solution of (5.1) tends exponentially to the periodic solution

$$x(t) = p(\omega,\varphi) \sin(\omega t) + q(\omega,\varphi) \cos(\omega t), \text{ for } t \geq 0 \tag{5.14}$$

where

$$p(\omega,\varphi) := \left(A^2 + \omega^2 I\right)^{-1} \left(\omega \sin(\varphi) I - \cos(\varphi) A\right) B$$

$$q(\omega,\varphi) := -\left(A^2 + \omega^2 I\right)^{-1} \left(\sin(\varphi) A + \omega \cos(\varphi) I\right) B \tag{5.15}$$

It follows from (5.14), (5.15) that the following equation holds for the output of system (5.1)

$$|y(t)|^2 = \frac{1}{2}\left(|Cq(\omega,\varphi)|^2 + |Cp(\omega,\varphi)|^2\right) + \frac{1}{2}\left(|Cq(\omega,\varphi)|^2 - |Cp(\omega,\varphi)|^2\right) \cos(2\omega t) + p'(\omega,\varphi) C' C q(\omega,\varphi) \sin(2\omega t) \tag{5.16}$$

and consequently we obtain the following equation from (5.16) for every initial condition:

$$\limsup_{t \to +\infty} \left(|y(t)|\right) = R(\omega,\varphi) \tag{5.17}$$

where

$$R(\omega,\varphi) :=$$
$$\frac{\sqrt{2}}{2}\sqrt{\left(|Cq(\omega,\varphi)|^2 + |Cp(\omega,\varphi)|^2\right) + \sqrt{\left(|Cq(\omega,\varphi)|^2 - |Cp(\omega,\varphi)|^2\right)^2 + 4\left(p'(\omega,\varphi)C'Cq(\omega,\varphi)\right)^2}} \quad (5.18)$$

Since the solution with initial condition $x_0 \in \Re^n$ corresponding to input given by (5.13) will have the same $\limsup_{t \to +\infty}(|y(t)|)$ as the solution with initial condition $\exp(Ar)x_0 + \int_0^r \exp(A(r-s))B\sin(\omega s + \varphi)ds$ corresponding to input $u(t) = \sin(\omega t + \varphi + \omega r)$ (the semigroup property), it follows that $R(\omega,\varphi) = R(\omega, \varphi + \omega r)$ for every $r \geq 0$. Consequently, the function $R(\omega,\varphi)$ is independent of $\varphi \in \Re$. Thus we get from (5.15), (5.17), (5.18) (with $\varphi = 0$) and (2.13):

$$\gamma = \gamma_{as} \geq \sup_{\omega > 0}(\Psi(\omega)) \quad (5.19)$$

where
$$\Psi(\omega) :=$$
$$\frac{\sqrt{2}}{2}\sqrt{\left(\omega^2|C\xi(\omega)|^2 + |CA\xi(\omega)|^2\right) + \sqrt{\left(\omega^2|C\xi(\omega)|^2 - |CA\xi(\omega)|^2\right)^2 + 4\omega^2\left(\xi'(\omega)A'C'C\xi(\omega)\right)^2}} \quad (5.19)$$

$$\xi(\omega) := \left(A^2 + \omega^2 I\right)^{-1} B \quad (5.20)$$

*5.4. Upper bound for the Minimum ISS-gain*

Let $C_1,...,C_n \in \Re^n$ be an orthonormal basis of $\Re^n$. Equation (5.7) implies that for all $u \in L^\infty(\Re_+;\Re)$ with $\sup_{t \geq 0}(|u(t)|) = 1$ the following inequalities hold:

$$\sup_{t \geq 0}\left(\left|\int_0^t C_i'\exp(A(t-s))Bu(s)ds\right|\right) \leq \int_0^{+\infty}|C_i'\exp(As)B|ds, \text{ for } i = 1,...,n \quad (5.21)$$

Since $C_1,...,C_n \in \Re^n$ is an orthonormal basis of $\Re^n$, it follows that $|x| = \sqrt{\sum_{i=1}^n (C_i'x)^2}$ for all $x \in \Re^n$. Therefore, we get from (5.21) for all $u \in L^\infty(\Re_+;\Re)$ with $\sup_{t \geq 0}(|u(t)|) = 1$:

$$\left|\int_0^t \exp(A(t-s))Bu(s)ds\right| = \sqrt{\sum_{i=1}^n \left(\int_0^t C_i'\exp(A(t-s))Bu(s)ds\right)^2}$$
$$\leq \sqrt{\sum_{i=1}^n \left(\int_0^{+\infty}|C_i'\exp(As)B|ds\right)^2} \quad (5.22)$$

It follows from (5.22) that the minimum ISS-gain of (5.1) satisfies the estimate:

$$\gamma \leq \inf\left\{\sqrt{\sum_{i=1}^{n}\left(\int_{0}^{+\infty}|C_i'\exp(As)B|ds\right)^2} : C_1,...,C_n \in \Re^n \text{ is an orthonormal basis of } \Re^n\right\} \quad (5.23)$$

Using formula (5.23) with the usual orthonormal basis of $\Re^n$ we are in a position to estimate that the minimum ISS-gain for (5.1) with $B=(b_1,...,b_n)' \in \Re^n$, $A=-diag(\lambda_1,...,\lambda_n)$, $\lambda_i > 0$ ($i=1,...,n$) satisfies $\gamma \leq \sqrt{\frac{b_1^2}{\lambda_1^2}+...+\frac{b_n^2}{\lambda_n^2}}$. This gives the exact value of the minimum ISS-gain: system (5.1) in this case satisfies assumption (H) with $P=I_n$ and formula (5.12) holds. On the other hand, the standard estimation of the ISS-gain gives $\gamma \leq \int_0^{+\infty}|\exp(As)B|ds = \int_0^{+\infty}\sqrt{b_1^2\exp(-2\lambda_1 s)+...+b_n^2\exp(-2\lambda_n s)}ds$, which is (in general) a conservative estimation of the minimum ISS-gain.

### 5.5. Upper bound for the Minimum IOS-gain

An upper bound of the minimum IOS-gain can be found by means of formulas (4.6), (4.8). When a $T$-periodic input is applied to system (5.1) then all solutions exponentially tend to the $T$-periodic solution

$$x(t) = -(\exp(AT)-I)^{-1}\int_0^T \exp(A(T-s))Bu(s+t)ds \quad (5.24)$$

It follows from (5.24) that the following estimate holds for the output response of the $T$-periodic solution:

$$\max_{t\in[0,T]}(|y(t)|) \leq \int_0^T |C(\exp(AT)-I)^{-1}\exp(As)B|ds \sup_{t\in[0,T]}(|u(t)|) \quad (5.25)$$

Since all solutions exponentially tend to the $T$-periodic solution, it follows from (5.25) that for every $T$-periodic input with $\sup_{t\in[0,T]}(|u(t)|)$, the estimate below holds:

$$\limsup_{t\to+\infty}(|y(t)|) \leq \int_0^T |C(\exp(AT)-I)^{-1}\exp(As)B|ds \quad (5.26)$$

Consequently, we obtain from (4.6), (4.8) and (5.26):

$$\gamma \leq \sup_{T>0}\left(\int_0^T |C(\exp(AT)-I)^{-1}\exp(As)B|ds\right) \quad (5.27)$$

Formula (5.27) should be compared with the standard estimate of the IOS-gain $\gamma \leq \int_0^{+\infty}|C\exp(As)B|ds$.

# 6. Concluding Remarks

The paper presented a fundamental relation between OAG's and IOS-gains for linear systems. For any ISS, strictly causal linear system the minimum OAG is equal to the minimum IOS-gain (Theorem 4.3). An extension of this result to the nonlinear case would involve the definition of minimum IOS-gains and OAG's for an ISS control system. This is not easy since in the nonlinear case, the effect of the input is combined in a nonlinear way with the effect of the initial conditions. The extension of Theorem 4.3 to the nonlinear case is an open problem.